\documentstyle{article}
\makeatletter\@addtoreset{equation}{section}\makeatother

\begin{document}
\begin{center} 
{\Large\bf ASTROPHYSICAL THERMONUCLEAR FUNCTIONS}\\ [1cm]
{\bf W.J. Anderson$^1$, H.J. Haubold$^2$, and A.M. Mathai$^1$}
\end{center}
\noindent
$^1$ Department of Mathematics and Statistics, McGill
University,
Montreal, P.Q., Canada H3A 2K6\par
\medskip
\noindent
$^2$Office for Outer Space, United Nations, New York, N.Y. 10017,
USA\par
\bigskip
\noindent
Received...
\clearpage
\noindent
Abstract.\hspace{0.2cm}Stars are gravitationally stabilized fusion
reactors changing their chemical composition while transforming
light atomic nuclei into heavy ones. The atomic nuclei are supposed
to be in thermal equilibrium with the ambient plasma. The majority
of reactions among nuclei leading to a nuclear transformation are
inhibited by the necessity for the charged participants to tunnel
through their mutual Coulomb barrier. As theoretical knowledge and
experimental verification of nuclear cross sections increases it
becomes possible to refine analytic representations for nuclear
reaction rates. Over the years various approaches have been made to
derive closed-form representations of thermonuclear reaction rates
(Critchfield 1972, Haubold and John 1978, Haubold, Mathai and
Anderson 1987).
They show that
the reaction rate contains the astrophysical
cross section factor and its derivatives which has to be
determined experimentally, and an integral part of the
thermonuclear reaction rate independent from experimental results
which can be treated by closed-form representation techniques in
terms of generalized hypergeometric functions. In this paper 
mathematical/statisti
cal techniques for deriving closed-form
representations of thermonuclear functions,  particularly the four
integrals 
\begin{eqnarray*}
I_1(z,\nu)& \stackrel{def}{=} & \int_0^\infty y^\nu
e^{-y}e^{-zy^{-\frac{1}{2}}}dy,\\
I_2(z,d,\nu)& \stackrel{def}{=} & \int_0^d y^\nu
e^{-y}e^{-zy^{-\frac{1}{2}}}dy,\\
I_3(z,t,\nu)& \stackrel{def}{=} & \int_0^\infty y^\nu
e^{-y}e^{-z(y+t)^{-\frac{1}
{2}}}dy,\\
I_4(z,\delta,b,\nu) & \stackrel{def}{=} &\int_0^\infty y^\nu
e^{-y}e^{-by^\delta}e^{-zy^{-\frac{1}{2}}}dy,
\end{eqnarray*}
will be
summarized and numerical results for them will  
be given. The separation of thermonuclear functions from
thermonuclear reaction rates is our preferred result. The purpose
of the paper is also to compare numerical
results for approximate and closed-form  
representations of thermonuclear functions. This paper completes
the work of Haubold, Mathai, and Anderson (1987).
\clearpage
\noindent
\section{Barrier penetration at astrophysical energies}
The majority of nuclear reactions of astrophysical interest are
inhibited by the necessity for the charged participants to tunnel
through their mutual Coulomb barrier. Nuclear processes such as
$\alpha$-decay and decay by emission of heavier nuclei are also
mediated by penetration through a static, one-dimensional Coulomb
potential barrier. Barrier penentration factors in nuclear reaction
rates take into account the exponential nature of the tail of the
nuclear potential. A great impact on the nature of the potential
has also the inclusion of the electron screening of the reacting
particles, which leads to potentials of the Yukawa type, which
exhibits a change of the height and width of the barrier compared
to the Coulomb type of potential (Fowler 1984).\par
In order to extrapolate measured nuclear cross sections $\sigma(E)$
down to astrophysical energies, the nuclear cross section factor
$S(E)$ is introduced by
\begin{equation}
\sigma(E)=\frac{S(E)}{E} exp\left\{ -
2\pi\eta\right\},\label{eq:1.1}
\end{equation}
where $\eta$ is the Sommerfeld parameter
\begin{equation}
\eta = \frac{Z_1Z_2e^2}{\hbar
v}=\left(\frac{\mu}{2}^{1/2}\frac{Z_1Z_2e^2}{\hbar
E^{1/2}}\right),\label{eq:1.2} 
\end{equation}
with Z the atomic charge and v [E] the asymptotic relative velocity
[kinetic energy] of the reacting nuclei (Fowler 1984). Thus, the
cross section is
given by the product of the cross section factor to be determined
experimentally, the square of the de Broglie wavelength due to
quantum mechanics $(\sim E^{-1}),$ and the barrier penetration
factor. The quantity exp$\left\{-2\pi \eta \right\}$ takes
exclusively s-wave transmission into account, describing
penetration to the origin through a pure Coulomb potential. Nuclear
reactions rates are extremely sensitive to the precise numerical
value in the argument of this exponential factor. The inclusion of
uncertainties in the shape of the nuclear potential and
contributions from non s-wave transmission, respectively, are very
important for deriving specific nuclear reaction rates but do not
change the overall energy dependence of the nuclear cross-section
given in (1.1). Actually, uncertainties in the shape of the nuclear
potential tail and contributions from non s-wave terms are only
important for heavy-ion reactions. In the following we are focusing
on reaction rates of the proton capture type, i.e. small value of
the reduced mass $\mu$ and small value of the atomic charge product
$Z_1Z_2e^2$. The main uncertainty in (1.1) lies in the variation of
the cross section factor S(E) with energy, which depends primarily
on the value chosen for the radius at which formation of a compound
nucleus between two interacting nuclei or nucleons occurs (Brown
and Jarmie 1990).\par
The separation of the barrier penetration factor in (1.1) is based
on
the solution of the Schr\"{o}dinger equation for the Coulomb wave
functions. Therefore the cross section $\sigma(E)$ in (1.1) can be
parametrized even more precisely by either expanding $S(E)$ into a
Taylor series about zero energy because of its slow energy
dependence,
\begin{equation}
S(E)=S(0)\left[1+\frac{S'(0)}{S(0)}E+\frac{1}{2}\frac{S"(0)}{S(0)
}E^2\right],\label{eq:1.3}
\end{equation}
where S(0) is the value of S(E) at zero energy, and S'(0) and S"(0)
are the first and second derivatives of S(E) with respect to energy
evaluated at E=0, respectively, or to elaborate on the action
integral $I(r_1, r_2)$ to include effects due to the shape of the
nuclear potential and non s-wave contributions, where $r_1$ and
$r_2$ are the inner turning point and outer turning point,
respectively, where the reacting particles tunnel from $r_1$ to
$r_2$ in a Coulomb plus nuclear field (Smith, Kawano, and Malaney
1993). Then the barrier penetration
factor in (1.1) can be expressed in terms of this action integral
as
\begin{equation}
T_0=exp\left\{ -2I(r_1, r_2)\right\},\label{eq:1.4}
\end{equation}
which simplifies for a Coulomb field and for $r_1=0$ to be
$I_c(0,r_2)=\pi \eta$, where $I_c(0,r_2)$ is the sharp-cutoff
Coulomb integral. If one takes into account non s-wave terms and
does not confine to the sharp-cutoff approximation of the Coulomb
integral in (1.4), the overall energy dependence of the nuclear
cross
section $\sigma(E)$ can be approximated by
\begin{equation}
\sigma(E)=\frac{S(E)}{E}\left\{ C_1
\frac{1}{E^{1/2}}+C_2E^{1/2}+C_3(C_4+E^{1/2})+\ldots\right\},
\label{eq:1.5}
\end{equation}
where the leading term containing $C_1E^{-1/2}$ corresponds to the
exponential term in (1.1); $C_2,C_3$ and $C_4$ are energy
independent
nuclear constants (Rowley and Merchant 1991).\par
Electron screening of reacting nuclei brings about a considerable
enhancement of nuclear reactions, particularly in high-Z matter.
The Coulomb potential is modified by the presence of a polarising
cloud of electrons surrounding the positive ions. The potential
seen by a reacting nucleus is found to be narrower than the Coulomb
potential and quantum-mechanical tunneling through the barrier
becomes easier. The barrier penetration factor in (1.1) taking into
account a screened potential can be written in terms of a screening
parameter t,
\begin{equation}
\sigma(E)=\frac{S(E)}{E}exp\left\{ -
2\pi(\frac{\mu}{2})^{1/2}\frac{Z_1Z_2e^2}{\hbar
(E+t)^{1/2}}\right\},\label{eq:1.6} 
\end{equation}
where $t=Z_1Z2e^2K$ and $K$ denotes the Debye-H\"{u}ckel length.
Screened nuclear reaction rates are extremely sensitive to the
precise numerical value of the argument of the exponential factor
in
(1.6). 
\section{Evolution towards the Maxwellian equilibrium distribution}
It is a major assumption in deriving nuclear reaction rates that
the reacting nuclei are supposed to be in thermal equilibrium with
the ambient plasma. This assumption can be justified by comparing
the characteristic time for significant energy exchanges by Coulomb
collisions with the characteristic time it takes the nuclear
reaction to produce the final nucleus. Generally the Coulomb
collision time is many orders of magnitude smaller than the time to
produce the final nucleus which is the natural condition that the
nuclei are in thermal equilibrium with the ambient plasma. Thus the
velocity distribution function of nuclei is Maxwell-Boltzmannian.
The state of the plasma at time t is described by the distribution
function $nf(v,t)$, where n is the constant particle number
density, $\vec{v}$ is the velocity variable, and $v=\mid
\vec{v}\mid$. Conservation of mass and energy imply that
\[\int d^3v f(v,t) = 1,\]
\begin{equation}
\int d^3v v^2 f(v,t) = \frac{3kT}{\mu},\label{eq:2.1}
\end{equation}
where T is the constant kinetic temperature and $\mu$ denotes the
mass. In a gravitationally stabilized stellar fusion reactor, as
$t\rightarrow \infty$, $f(v,t)$ tends to the Maxwell-Boltzmann
distribution function, 
\begin{equation}
f(v,\infty) dv=\left(\frac{\mu}{2\pi kT}\right)^{3/2} exp\left\{
-\frac{\mu v^2}{2kT}\right\} 4\pi v^2 dv.\label{eq:2.2}
\end{equation}
In a thermonuclear plasma, the reaction rate arises from an
integral of the nuclear cross section (equations (1.1) or (1.6)),
times
velocity, times the Maxwell-Boltzmann distribution of velocities
(2.2), 
\begin{equation}
<\sigma v>=\left(\frac{\mu}{2\pi kT}\right)^{3/2}
\int_0^\infty dv \sigma (v) v^3 exp\left\{ -\frac{\mu
v^2}{2kT}\right\}.\label{eq:2.3}
\end{equation}
It is evident from (2.3) that the kernel of the integral consists
of
a product of the steeply falling Maxwell-Boltzmann distribution
(2.2) and
the rapidly rising cross section (1.1) or (1.6) to produce a not
quite
symmetrical peak, commonly called the Gamow peak. This peak
justifies the fact, that reaction rates are extremely sensitive to
the precise numerical values in the arguments of the exponential
factors exhibiting the exponential nature of the tail of the
nuclear potential and the exponential nature of the tail of the
velocity distribution function.\par
The Maxwell-Boltzmann distribution is a solution of the general
nonlinear Boltzmann equation which itself reveals as notoriously
complicated. The system of particles here is considered to be an
infinite, spatially homogeneous and isotropic gas containing a
variety of nuclei. It is also assumed that only binary reactions
need to be taken into account, so that the Boltzmann equation
applies. Additionally the assumption is made that the nuclear
reactions are isotropic, i.e., the cross section $\sigma$ is
independent of the collision angle. Maxwell established that the
low-order moments of the distribution function effectively relax
toward their equilibrium values in just a few mean  collision
times. This corresponds, as discussed before, to the property that
the low-energy part of the distribution attains Maxwell-
Boltzmannian form in such a time interval. Nonlinear relaxation has
been discussed by Kac (1955).\par
On several occasions the question has been raised whether there may
exist intermediate distributions that will evolve in such a way
that the high-velocity tail of the respective velocity
distribution will, at certain typically high velocities and for
certain time-intervals, display significant enhancement or
depletion with respect to the steady-state Maxwell-Boltzmann
distribution. Such a modification of the tail away from the
Maxwell-Boltzmann distribution would significantly change the Gamow
peak in (2.3) and subsequently would alter the respective reaction
rates among nuclei in the plasma of the gravitationally stabilized
stellar fusion reactor. A certain type of nonequilibrium
distribution functions have been studied by Krook and Wu (1976,
1977), Tjon and Wu (1979), and Barnsley and Cornille (1981) by
investigating solutions of the Boltzmann equation which approach an
equilibrium distribution when $t\rightarrow \infty$ in a nonuniform
fashion. This nonuniformity is due to the high velocity tail of the
distribution and indicates that linearization techniques can not be
fully justified for high velocities even when the state of the
physical
system is close to Maxwell-Boltzmannian behavior. Their model
considerations, while studying the relaxation of solutions of the
Boltzmann equation towards the steady-state Maxwell-Boltzmann
distribution, encourage the investigation of reaction rates
containing a modified Maxwell-Boltzmann distribution.
Having discussed the energy dependence of the nuclear cross section
in Section 1 and Maxwell-Boltzmann distribution in Section 2,
respectively, the following four integrals can be derived,
representing thermonuclear functions for four quite different
physical conditions. The standard case of the thermonuclear
function contains the nuclear cross section (1.1), the energy
dependent term of the Taylor series in (1.3), and the steady-state
Maxwell-Boltzmann distribution function (2.2),
\begin{equation}
I_1(z,\nu) \stackrel{def}{=} \int_0^\infty y^\nu
e^{-y}e^{-zy^{-\frac{1}{2}}}dy\label{eq:2.4}
\end{equation}
where $y=E/kT$ and $z=2\pi(\mu/2kT)^{1/2}Z_1Z_2e^2/\hbar$.
Considering dissipative collision processes in the thermonuclear
plasma cut off of the high energy tail of the Maxwell-Boltzmann
distribution may occur, thus we write for (2.4),
\begin{equation}
I_2(z,d,\nu) \stackrel{def}{=}  \int_0^d y^\nu
e^{-y}e^{-zy^{-\frac{1}{2}}}dy,\label{eq:2.5}
\end{equation}
where d denotes a certain typically high energy.\par
Accomodating screening effects in the standard thermonuclear
function we have to use the nuclear cross section (1.6) and the
steady-state Maxwell-Boltzmann distribution function (2.2) which
leads to 
\begin{equation}
I_3(z,t,\nu) \stackrel{def}{=} \int_0^\infty y^\nu
e^{-y}e^{-z(y+t)^{-\frac{1}
{2}}}dy,\label{eq:2.6}
\end{equation}
where t is the electron screening parameter. \par
Finally, if due to plasma effects a depletion of the Maxwell-
Boltzmann distribution has to be taken into account, the
thermonuclear function can be written in the follwing form
\begin{equation}
I_4(z,\delta,b,\nu) \stackrel{def}{=} \int_0^\infty y^\nu
e^{-y}e^{-by^\delta}e^{-zy^{-\frac{1}{2}}}dy, \label{eq:2.7}    
\end{equation}
where the parameter $\delta$ exhibits the enhancement or reduction
of the high-energy tail of the Maxwell-Boltzmann distribution.\par
In the following Sections mathematical/statistical techniques for
deriving closed-form representations of the four thermonuclear
functions (2.4) - (2.6) will be summarized, their asymptotic forms
will be given and numerical results for both of them derived.
\section {Mathematical preliminaries}
First of all, we need to recall the {\sl gamma function}, defined
for complex $z$ by 
\[\Gamma(z)=\int_0^\infty t^{z-1}e^{-t}dt, \Re(z)>0.\]
The definition of $\Gamma(z)$ can be extended to the entire complex
plane where it is analytic
except for simple poles at $0$ and the negative real integers. An
important property we shall need
is the {\sl multiplication formula}
\begin{equation}
\Gamma(mz)=(2\pi)^{\frac{1-m}{2}}m^{mz-\frac{1}
{2}}\Gamma(z)\Gamma(z+\frac{1}
{m})\cdots\Gamma(z+\frac{m-1}{m}),\label{eq:3.1}
\end{equation}
which is valid for all $z$ and all integers $m\ge1$.\par
\noindent
{\bf Definition.} The function
\begin{equation}
G^{m,n}_{p,q}(z)=G^{m,n}_{p,q}\left(z\bigg|^{a_1,\ldots,a_p}_{b_1
,\ldots b_q}\right)=
\frac{1}{2\pi
i}\int_L\frac{\Pi_{j=1}^m\Gamma(b_j+s)\Pi_{j=1}^n\Gamma(1-a_j-s)}
{
\Pi_{j=m+1}^q\Gamma(1-b_j-s)\Pi_{j=n+1}^p\Gamma(a_j+s)}z^{-s}ds,
z\neq 0,\label{eq:3.2}
\end{equation}
is called the {\sl $G$-function} and is originally due to Meijer
(cp. Mathai and Saxena 1973).
Here, $i=\sqrt{-1}$;
$m$, $n$, $p$, and $q$ are integers with $0\leq n \leq p$ and
$0\leq
m\leq q$. In (3.2), and
throughout this paper, an empty product is interpreted as unity
(similarly an empty sum as zero). The 
$a_j,j=1,\ldots,p$ and $b_j,j=1,\ldots,q$ are complex numbers such
that no pole of $\Gamma(b_j+s),
j=1,\ldots,m$ coincides with any pole of $\Gamma(1-a_j-s),
j=1,\ldots,n$.  $L$ is a contour
separating the poles of $\Gamma(b_j+s), j=1,\ldots,m$ from the
poles of $\Gamma(1-a_j-s),
j=1,\ldots,n$. 
At this point, it is not clear that the integral in (3.2)
even exists. Conditions
on the contour $L$ and the various parameters must be imposed in
order that the integral
converges. These conditions, as well as properties of the
$G$-function may be found in
Luke (1969), chapter 5. However, for the $G$-functions encountered
in this paper, it suffices to
know that the integral in (3.2) is well-defined for all
$z\neq 0$ if  

\begin{enumerate}
\item[(i)] $L$ is a loop beginning and ending at $-\infty$ and
encircling all poles of
$\Gamma(b_j+s), j=1,\ldots,m$, once in the positive direction, but
none of the poles of
$\Gamma(1-a_j-s), j=1,\ldots,n$, and 
\item[(ii)] $q\geq 1$ and $p<q$.
\end{enumerate}
Moreover, under these conditions the integral can be evaluated as
a sum of residues at the poles
of $\Gamma(b_j+s), j=1,\ldots,m$.

One property that we will certainly require in the sequel is the
asymptotic behaviour of $G^{q,0}_{p,q}(z)$ as $|z|\rightarrow
\infty$. From Luke (1969)
page 179, we have
\begin{equation}
G^{q,0}_{p,q}\left(z\bigg|^{a_1,\ldots,a_p}_{b_1,\ldots
b_q}\right)\sim\frac{(2\pi)^{\frac{\sigma-1}
{2}}}{\sigma^{\frac{1}{2}}}e^{-\sigma
z^{\frac{1}{\sigma}}}z^\theta \; \mbox{as}\;|z|\rightarrow
\infty,\;
|\arg z|\leq(\sigma+\epsilon)\pi-\frac{\sigma}{2},\label{eq:3.3}
\end{equation}
where 
\[\sigma=q-p>0,\epsilon=\left\{ \begin{array}{lll} 
\frac{1}{2} & \mbox{if $\sigma=1$,}\\
1 & \mbox{if $\sigma\geq
1$,}
\end{array}
\mbox{and}\;\sigma\theta=\frac{1
}
{2}(1-\sigma)+\sum_{j=1}^q
b_j-\sum_{j=1}^p a_j.\right.\]
\medskip
\noindent
{\bf Definition.} Let $f(t)$ be a function defined for $t>0$.
Then
\begin{eqnarray}
M_f(s)&=&\int_0^\infty t^{s-1}f(t)dt,
\alpha<\Re(s)<\beta,\label{eq:3.4}\\
f(t)&=&\frac{1}{2\pi
i}\int_{c-i\infty}^{c+i\infty}t^{-s}M_f(s)ds,\label{eq:3.5}
\end{eqnarray}
is called a {\it Mellin transform pair}. (3.4) is called
the {\it Mellin transform}, and
(3.5) is the inversion formula. The transform normally
exists only in the strip
$\alpha<\Re(s)<\beta$, and the inversion contour must lie in this
strip.\par
\vspace{0.3cm}
\noindent
{\bf Lemma 3.1} Let $f_1(t)$ and $f_2(t)$ be two functions with
{\it Mellin transforms} $M_{f_1}(s)$ and
$M_{f_2}(s)$. Then
\begin{equation}
\int_0^\infty v^{-1}f_1(v)f_2(\frac{u}{v})\,dv=\frac{1}{2\pi
i}\int_L
M_{f_1}(s)M_{f_2}(s)u^{-s}ds.
\label{eq:3.6}
\end{equation}
{\bf Proof.} Our proof is statistical. We suppose that
$f_1(t)\ge0$,
$f_2(t)\ge0$, $\int_0^\infty
f_1(t)dt<\infty$, and $\int_0^\infty f_2(t)dt<\infty$ (the
application below will satisfy these
criteria). By scaling if necessary, we can assume that $f_1(t)$ and
$f_2(t)$ are density functions.
Let $X$ and $Y$ be independent random variables having density
functions $f_1(t)$ and $f_2(t)$
respectively. Then the left-hand side of (3.6) is the
density function $g(u)$ of the random
variable $U=XY$. Let us look at the right-hand side. We have
$M_{f_1}(s)=E(X^{s-1})$ and
$M_{f_2}(s)=E(Y^{s-1})$, and therefore 
$$M_g(s)=E(U^{s-1})=E(X^{s-1}Y^{s-1})=E(X^{s-1})E(Y^{s-1})=M_{f_1
}(s)M_{f_2}(s).$$
It follows that the right-hand side of (3.6) is ${1\over
2\pi i}\int_L M_g(s)u^{-s}ds$,
which is the formula for the inverse Mellin transform of $M_g(s)$.
Thus the right-hand side is
$g(u)$ as well.
\section{Representation of the four integrals in terms of
$G$-functions}
{\bf Theorem 4.1 (Saxena (1960), Mathai and Haubold (1988))} For
$z>0,
p>0, \rho\leq 0$, and integers
$m,n\geq1$, we have 
\begin{eqnarray}
p\int_0^\infty t^{-n\rho}e^{-pt}e^{-zt^{-\frac{n}
{m}}}dt & = & p^{n\rho}(2\pi)^{\frac{1}{2}(2-n-
m)}m^{\frac{1}{2}}n^{\frac{1}{2}-n\rho} \nonumber \\
& \times & G_{0,m+n}^{m+n,0}\left(\frac{z^mp^n}
{m^mn^n}\bigg|_{0,\frac{1}{m},\ldots,\frac{m-1}{m};\frac{1-n\rho}
{n},\ldots,\frac{n-n\rho}{n}}\right)\label{eq:4.1}
\end{eqnarray}
{\bf Proof.} Define $f_1(t)=t^{1-n\rho}e^{-t}$ and
$f_2(t)=e^{-t^{\frac{n}{m}}}$ for $t>0$. Then the Mellin
transforms are
\[M_{f_1}(s)=\int_0^\infty t^{s-1}f_1(t)dt=\int_0^\infty
t^{1-n\rho+s-1}e^{-t}dt=\Gamma(1-n\rho+s), \Re(1-n\rho+s)>0,\]
and 
\[M_{f_2}(s)=\int_0^\infty t^{s-1}f_2(t)dt=\int_0^\infty
t^{s-1}e^{-t^{\frac{n}{m}}}dt=\frac{m}
{n}\Gamma(\frac{m}{n}s),\Re(s)>0.\]
Then by setting $v=pt$ and $u=z^{\frac{m}{n}}p$, and using the
lemma,
we have
\begin{eqnarray}
p\int_0^\infty t^{-n\rho}e^{-pt}e^{-zt^{-\frac{n}{m}}}dt & =
& p^{n\rho}\int_0^\infty 
v^{-n\rho}e^{-v}e^{-(\frac{u}{v})^\frac{n}{m}}dv=p^{n\rho}\int_0^
\infty
v^{-1}f_1(v)f_2(\frac{u}{v})dv \nonumber \\
& = & \frac{p^{n\rho}}{2\pi i}\int_L
M_{f_1}(s)M_{f_2}(s)u^{-s}ds=\frac{p^{n\rho}}{2\pi i}\int_L
\Gamma(1-n\rho+s)\frac{m}{n}\Gamma(\frac{m}{n}s)u^{-s}ds
\nonumber \\
& = & \frac{mp^{n\rho}}{2\pi i}\int_{L'}
\Gamma(1-n\rho+ns')\Gamma(ms')(z^mp^n)^{-s'}ds,\label{eq:4.2}
\end{eqnarray}
where we made a change of variable $s=ns'$. The $G$-function
appearing on the right-hand side of
(4.1) is 
\begin{eqnarray}
& G_{0,m+n}^{m+n,0}\left(\frac{z^mp^n}
{m^mn^n}\bigg|_{0,\frac{1}{m},\ldots,\frac{m-1}{m};\frac{1-n\rho}{
n},\ldots,\frac{n-n\rho}{
n}}\right) \nonumber \\
& = \frac{1}{2\pi i}\int_L\Gamma(s)\Gamma(\frac{1}
{m}+s)\cdots\Gamma(\frac{m-1}
{m}+s)\Gamma(\frac{1-n\rho}{n}+s)\cdots\Gamma(\frac{n-n\rho}
{n}+s)\left(\frac{z^mp^n}{m^mn^n}\right)^{-s}ds.\label{eq:4.3}
\end{eqnarray}
By the multiplication formula in (3.1), we have
\begin{equation}
\Gamma(1-n\rho+ns)=\Gamma(n[\frac{1}{n}-\rho+s])=(2\pi)^{\frac{1-
n}{2}}n^{n(\frac{1}{n}-\rho+s)-\frac{1}{
2}}\Gamma(\frac{1-n\rho}{n}+s)\cdots\Gamma(\frac{n-n\rho}
{n}+s).\label{eq:4.4}
\end{equation}
Thus by applying the multiplication formula and (4.4) to
(4.3),
we get 
\begin{eqnarray}
G_{0,m+n}^{m+n,0}\left(\frac{z^mp^n}
{m^mn^n}\bigg|_{0,\frac{1}{m},\ldots,\frac{m-1}{m};\frac{1-n\rho}
{n},\ldots,\frac{n-n\rho}
{n}}\right) & = & \frac{(2\pi)^{\frac{m+n}{2}-1}m^\frac{1}
{2}n^{n\rho-\frac{1}
{2}}}{2\pi i} \nonumber \\
&\times &
\int_L\Gamma(ms)\Gamma(1-n\rho+ns)(z^mp^n)^{-s}ds.\nonumber \\
& & \hspace{5.5cm}\label{eq:4.5}
\end{eqnarray}
By comparing (4.2) and (4.5), we obtain
(4.1).
\vspace{.2in}
By setting $m=2,n=1,p=1,$ and $\rho=-\nu$, we obtain\par
\noindent
{\bf Corollary 4.2} For $z>0$ and $\nu\geq 0$, we have
\begin{equation}
I_1(z,\nu)=\int_0^\infty y^\nu e^{-y}e^{-zy^{\frac{1}
{2}}}dy=\pi^{\frac{-1}
{2}}G^{3,0}_{0,3}\left(\frac{z^2}{4}\bigg|_{0,\frac{1}
{2},1+\nu}\right).\label{eq:4.6}
\end{equation}
\bigskip\noindent
The proof of the following theorem is similar to that of theorem
4.1.\par
\noindent
{\bf Theorem 4.3 (Mathai and Haubold (1988))} For $z>0, d>0, a>0$,
and
integers $m,n\geq 1$, we have
\begin{eqnarray*}
\int_0^dt^{-n\rho}e^{-at}e^{-zt^{-\frac{n}{m}}}dt &
= & \frac{m^{\frac{1}{
2}}}{n}(2\pi)^{\frac{1-m}
{2}}d^{1-n\rho}\\
& \times & \sum_{r=0}^\infty \frac{(-ad)^r}
{r!}G^{m+n,0}_{n,m+n}\left(\frac{z^m}
{d^nm^m}\bigg|^{-\rho+\frac{r+2}{n}+\frac{j-1}{n},j=1,\ldots,n}_
{-\rho+\frac{r+1}{n}+\frac{j-1}{n},j=1,\ldots,n;\frac{j-1}
{m},j=1,\ldots,m}\right)
\end{eqnarray*}
\vspace{.2in}
By setting $m=2,n=1,a=1,$ and $\rho=-\nu$, we obtain\par
\noindent
{\bf Corollary 4.4} For $z>0,d>0$, and $\nu\geq0$, we have 
\begin{equation}
I_2(z,d,\nu)=\int_0^d y^\nu e^{-y}e^{-zy^{\frac{-1}
{2}}}dy=\frac{d^{1+\nu}}{\pi^{\frac{1}
{2}}}\sum_{r=0}^\infty\frac{(-d)^r}{r!}G^{3,0}_{1,3}\left(\frac{z
^2}
{4d}\bigg|_{\nu+r+1,0,\frac{1}{2}}^{\nu+r+2}\right).\label{eq:4.7}
\end{equation} 
\bigskip\noindent
The integral $I_3$ may be worked out in terms of $I_1$ and $I_2$ as
follows. We have
\begin{eqnarray}
I_3(z,t,\nu) & \stackrel{def}{=} &\int_0^\infty y^\nu
e^{-y}e^{-z(y+t)^{\frac{-1}
{2}}}dy=\int_t^\infty(u-t)^\nu e^{-(u-t)}e^{-zu^{\frac{-1}{2}}}du
\mbox{(where $u=y+t$)}\nonumber\\
& = & e^t\int_t^\infty\sum_{r=0}^\nu(^\nu
_r)u^r(-t)^{\nu-r}e^{-u}e^{-zu^{\frac{-1}{2}}}du=
e^t\sum_{r=0}^\nu(^\nu_ r)(-t)^{\nu-r}\int_t^\infty
u^re^{-u}e^{-zu^{\frac{-1}{2}}}du\nonumber \\
& = & e^t\sum_{r=0}^\nu(^\nu
_r)(-t)^{\nu-r}\left[I_1(z,r)-I_2(z,t,r)\right]\label{eq:4.8}
\end{eqnarray}
\noindent
and
\begin{eqnarray}
I_4(z,\delta,b,\nu)& \stackrel{ref}{=} & \int_0^\infty y^\nu
e^{-y}e^{-by^\delta}e^{-zy^\frac{-1}{2}}dy=
\int_0^\infty y^\nu e^{-y}\sum_{r=0}^\infty \frac{(-b)^r}
{r!}y^{r\delta}e^{-zy^{\frac{-1}
{2}}}dy\nonumber \\
& = & \sum_{r=0}^\infty\frac{(-b)^r}{r!}\int_0^\infty
y^{\nu+r\delta}e^{-y}
e^{-zy^{\frac{-1}{2}}}dy=\sum_{r=0}^\infty\frac{(-b)^r}
{r!}I_1(z,\nu+r\delta).\label{eq:4.9}
\end{eqnarray}
In order to confidently exchange the summation and integral signs
in (4.9), the
quantity 

\[\int_0^\infty y^\nu e^{-y}e^{by^\delta}e^{-zy^{-{1\over 2}}}dy=
\int_0^\infty y^\nu e^{-y}\sum_{r=0}^\infty\left|{(-b)^r\over
r!}y^{r\delta}e^{-zy^{-{1\over
2}}}\right|dy\]
must be finite (by Fubini's theorem). Hence we expect the expansion
in (4.9) may not
be valid for large $b$ and $\delta$. (This was in fact borne out by
later numerical computations.)

To end this section, we use (3.3) to obtain asymptotic
formulas for the four integrals. By
a direct application of (3.3) to (4.6),
(4.7),
(4.8), and (4.9) (and some algebra in the case
of the last three), we obtain 
\begin{eqnarray}
I_1(z,\nu) & \sim & 2\left({\pi\over 3}\right)^{1\over
2}\left({z^2\over 4}\right)^{2\nu+1\over
6}e^{-3\left({z^2\over 4}\right)^{1/3}},\label{eq:4.10}\\
I_2(z,d,\nu) & \sim & d^{\nu+1}e^{-d}\left({z^2\over
4d}\right)^{-1/2}e^{-2\left({z^2\over
4d}\right)^{1/2}},\label{eq:4.11}\\
I_3(z,t,\nu) & \sim & 2\left({\pi\over 3}\right)^{1\over
2}e^t\left({z^2\over 4}\right)^{1\over 6}e^{-3 
\left({z^2\over 4}\right)^{1/3}}\left[\left({z^2\over
4}\right)^{1\over
3}-t\right]^\nu,\label{eq:4.12}\\
I_4(z,\delta,b,\nu) &\sim & 2\left({\pi\over 3}\right)^{1\over
2}\left({z^2\over 4}\right)^{2\nu+1\over
6}e^{-3\left({z^2\over 4}\right)^{1/3}}e^{-b\left({z^2\over
4}\right)^{\delta/3}},\label{eq:4.13}
\end{eqnarray}

all as $z\to\infty$.

\section{Series representations for the four integrals}
Series expressions for the four integrals can now be obtained by
evaluating the $G$-functions
using residue calculus. We will illustrate the method by doing this
for the integral $I_1(z,\nu)$.
This means that we have to evaluate the complex integral
\begin{equation}
G^{3,0}_{0,3}\left({z^2\over 4}\bigg|_{0,{1\over
2},1+\nu}\right)={1\over 2\pi i}\int_L 
\Gamma(s)\Gamma({1\over 2}+s)\Gamma(1+\nu+s)\left({z^2\over
4}\right)^{-s}ds.\label{eq:5.1}
\end{equation}
As previously mentioned, the right-hand side will be the sum
($R_1+R_2+R_3$ below) of the residues
of the integrand. We will assume that $\nu$ is a non-negative
integer (the analysis is slightly
different otherwise, as seen in Mathai and Haubold (1988)).  Then
the poles of the gammas in the
integrand of (5.1) are as follows: 
\begin{description}
\item[Poles of $\Gamma(s)$:] $s=0,-1,-2,\ldots$
\item[Poles of $\Gamma({1\over 2}+s)$:] $s=-{1\over 2},-{3\over
2},-{5\over 2},\ldots$
\item[Poles of $\Gamma(1+\nu+s)$;] $-\nu-1,-\nu-2,-\nu-3,\ldots$.
\end{description}
Note that $\Gamma(s)$ and $\Gamma(1+\nu+s)$ have some poles in
common. These will be poles of
order two. Thus the poles 
\begin{eqnarray*}
s & = & 0,-1,-2,\ldots,-\nu \mbox{are of order $1$ each,}\\
s & = & -{1\over 2},-{3\over 2},-{5\over 2},\ldots\mbox{are of
order $1$ each,}\\
s & = & -\nu-1,-\nu-2,-\nu-3,\ldots\mbox{are of order $2$ each.} 
\end{eqnarray*}
Using the facts that
\[lim_{s\to -r}(s+r)\Gamma(s)={(-1)^r\over r!},\;
\Gamma(a-r)={(-1)^r\Gamma(a)\over(1-a)_r},
\; r=0,1,2,\ldots;\; \Gamma\left({1\over
2}\right)=\pi^{1\over 2},\]

where 
\[(a)_r=\left\{ \begin{array}{ll} a(a+1)\cdots(a+r-1) & \mbox{if
$r\ge1$,}\\
1 & \mbox{if $r=0$,}
\end{array}
\right.\]
we find that the sum of residues of the integrand at the poles
$s=0,-1,\ldots,-\nu$ is
\begin{eqnarray*}
R_1 & = & \sum_{r=0}^\nu \lim_{s\to -r}(s+r)\Gamma(s)\Gamma({1\over
2}+s)\Gamma(1+\nu+s)
\left({z^2\over 4}\right)^{-s}\\
& = & \sum_{r=0}^\nu {(-1)^r\over r!}\Gamma({1\over
2}-r)\Gamma(1+\nu-r)
\left({z^2\over 4}\right)^r\\
& = & \pi^{1\over 2}\Gamma(1+\nu)\sum_{r=0}^\nu
{1\over\left({1\over
2}\right)_r(-\nu)_rr!}
\left(-{z^2\over 4}\right)^r.
\end{eqnarray*}
In exactly the same way, the sum of the residues at the poles
$s=-{1\over 2},-{3\over 2},-{5\over
2},\ldots$ is
\[R_2=-2\pi^{1\over 2}\Gamma\left({1\over
2}+\nu\right)\left({z^2\over 4}\right)^{1\over
2}{_0F_2}\left(-;{3\over 2},{1\over 2}-\nu;-{z^2\over 4}\right),\]
where ${_0F_2}$ is the hypergeometric function defined by
\[{_0F_2}(-;a,b;x)=\sum_{r=0}^\infty {x^r\over (a)_r(b)_rr!}.\]
Finally, the sum of the residues at
the poles $s=-\nu-1,-\nu-2,\ldots$ (each of order $2$) is
\begin{eqnarray*}
R_3&=\displaystyle\sum_{r=0}^\infty\lim_{s\to
-\nu-1-r}{\partial\over\partial
s}\left[(s+1+\nu+r)^2
\Gamma(s)\Gamma({1\over 2}+s)\Gamma(1+\nu+s)\left({z^2\over
4}\right)^{-s}\right]\\
&=\left({z^2\over
4}\right)^{1+\nu}\displaystyle\sum_{r=0}^\infty\left({z^2\over
4}\right)^r
\left[-\log\left({z^2\over 4}\right)+A_r\right]B_r,
\end{eqnarray*}
where
\begin{eqnarray}
A_r &=\left[1+{1\over 2}+\cdots+{1\over r}\right]+\left[1+{1\over
2}+\cdots+{1\over
r+\nu+1}\right]\\ \nonumber
& \hspace{0.5in}+\left[{1\over 1/2}+{1\over
3/2}+\cdots+{1\over (1/2)+\nu+r}\right]-3\gamma-2\log
2\label{eq:5.2}
\end{eqnarray}
(we take $[1+{1\over 2}+\cdots+{1\over r}]$ to be zero if $r=0$),
$\gamma=0.5772156649...$ is Euler's constant, and  

\begin{equation}
B_r={(-1)^{1+\nu+r}\Gamma\left(-{1\over 2}-\nu\right)\over
r!(r+\nu+1)!\left({3\over 2}+\nu\right)_r}.\label{eq:5.3} 
\end{equation}
By summing $R_1, R_2$, and $R_3$ and using (4.6) and
(5.1), we obtain the
following theorem.\\ 
{\bf Theorem  5.1} Let $\nu\ge0$ be an integer. Then for $z>0$, we
have
\begin{eqnarray*}
I_1(z,\nu)&=\Gamma(1+\nu)\displaystyle\sum_{r=0}^\nu{1\over\left(
{1\over
2}\right)_r(-\nu)_rr!}
\left(-{z^2\over 4}\right)^r-2\Gamma\left({1\over
2}+\nu\right)\left({z^2\over 4}\right)^{1\over
2}{_0F_2}\left(-;{3\over 2},{1\over 2}-\nu;-{z^2\over 4}\right)\\
&\hspace{.25in}+\pi^{-{1\over 2}}\left({z^2\over
4}\right)^{1+\nu}\displaystyle\sum_{r=0}^\infty\left({z^2\over
4}\right)^r \left[-\log\left({z^2\over 4}\right)+A_r\right]B_r
\end{eqnarray*}
where $A_r$ and $B_r$ are given above in (5.2) and
(5.3).
\bigskip\noindent
The integral $I_2(z,d,\nu)$ can be treated in the same way with the
following result.\par
\noindent
{\bf Theorem 5.2} Let $z>0$, $t>0$, and let $\nu\ge0$ be an
integer.
Then
\begin{eqnarray*}
I_2(z,d,\nu)&=&d^{\nu+1}\displaystyle\sum_{r=0}^\infty{(-d)^r\over
r!}\bigg\{-2\left({z^2\over 4d}\right)^{1\over
2}\displaystyle\sum_{\ell=0}^\infty{\left({z^2\over
4d}\right)^\ell\over
\ell!\left({3\over
2}\right)_\ell(\nu+r+{1\over 2}-\ell)}\\
&+&\displaystyle\sum_{\ell=0\atop\ell\neq\nu+r+1}^\infty
{\left({z^2\over 4d}\right)^\ell\over \ell!\left({1\over
2}\right)_\ell(\nu+r+1-\ell)}
+{2\left({z^2\over 4d}\right)^{\nu+r+1}\over(\nu+r+1)!\left({3\over
2}\right)_{\nu+r}}\left[-\log\left({z^2\over
4d}\right)+A\right]\bigg\}
\end{eqnarray*}
where
\[A=-2\gamma-2\log 2+\left[1+{1\over 2}+\cdots+{1\over
\nu+r+1}\right]+\left[{1\over 1/2}+{1\over
3/2}+\cdots+{1\over(1/2)+\nu+r}\right].\]

Since the integrals $I_3$ and $I_4$ have been expressed in
(4.8) and (4.9) in
terms of $I_1$ and $I_2$, similar expansions can and have been
derived for $I_3$ and $I_4$. However,
the exact details will not be given here.
\section {Computations and Conclusion}
Numerical computations for the series
expansions obtained above of the four
integrals $I_1$, $I_2$, $I_3$, and $I_4$ were made and compared to
the corresponding
approximations for large $z$ in (4.10)-(4.13).
The programming was carried out in
Pascal on a Macintosh II computer with a  numerical coprocessor,
for a wide range of parameter
values. Some of the results are shown in figures 1 to 4. Obviously,
the series computations fail
for large values of $z$. Since much effort was made in optimizing
the program for accuracy and  
countering problems of underflow and overflow, it is thought that
this failure is a result of machine and  
compiler numerical accuracy. It is evident, however, that the
missing portions of the \lq\lq exact\rq\rq  
curves can be replaced by the \lq\lq approximate\rq\rq curves.

For the sake of comparison, the four integrals were computed for
the same parameter values using the
numerical integration routines in Mathematica (Wolfram 1991). The
results were
identical to the results of the
previous paragraph, except that computations for larger values of
$z$ were possible. The results,
together with corresponding approximations, are plotted in figures
5 to 8.\par
\vspace{2cm}
\begin{center}
Acknowledgement
\end{center}  
The authors would like to thank the Natural Sciences and the
Engineering Research Councel of Canada for financial assistance for
this research project.
\clearpage
\noindent
References\par
\bigskip
\noindent
Barnsley, M., Cornille, H.: 1981, Proc. R. Soc. Lond. {\bf A374},
371\par
\medskip
\noindent
Brown, R.E., Jarmie, N.: 1990, Phys. Rev. {\bf C41}, 1391\par
\medskip
\noindent
Critchfield, C.L.: 1972, Analytic forms of the thermonuclear
function.\par 
In: {\sl Cosmology, Fusion, and Other  
Matters. George Gamow
Memorial\par 
Volume}, Edited by F. Reines, University of Colorado
Press, Colorado,\par 
pp. 186-191\par
\medskip
\noindent
Fowler, W.A.: 1984, Rev. Mod. Phys. {\bf 56}, 149\par
\medskip
\noindent
Haubold, H.J., John, R.W.: 1978, Astron. Nachr. {\bf 299}, 225\par
\medskip
\noindent
Haubold, H. J., Mathai, A. M., Anderson, W. J.: 1987, 
Thermonuclear\par 
functions. In: Proceedings
of the Workshop on Nuclear Astrophysics,\par 
Edited by W. Hillebrandt, R. Kuhfuss, E. Mueller,
J.W. Truran,\par 
{\sl Lecture Notes in
Physics Vol. 287}, Springer-Verlag,
Berlin pp. 102-110\par
\medskip
\noindent
Kac, M.: 1955, Foundations of kinetic theory. In: Proceedings\par
of the Third Berkeley Symposium on Mathematical Statistics and\par
Probability, University of California Press, Berkeley, pp. 171-
197\par

\medskip
\noindent
Krook, M., Wu, T.T.: 1976, Phys. Rev. Lett. {\bf 36}, 1107\par
\medskip
\noindent
Krook, M., Wu, T.T.: 1977, Phys. Fluids {\bf 20}, 1589\par
\medskip
\noindent
Luke, Y. L.: 1969, {\sl The Special Functions and Their
Approximations,\par 
Volume I}, Academic Press,
New York\par
\medskip
\noindent
Mathai, A. M., Haubold, H.J.: 1988, {\sl Modern Problems in
Nuclear and\par 
Neutrino
Astrophysics}, Akademie-Verlag, Berlin\par
\clearpage
\noindent
Mathai, A.M., Saxena, R.K.: 1973, Generalized Hypergeometric\par
Functions with Applications in Statistics and Physical
Sciences,\par
Lecture Notes in Mathematics Vol. 348, Springer-Verlag, Berlin\par
\medskip
\noindent
Rowley, N., Merchant, A.C.: 1991, Astrophys. J. {\bf 381}, 591\par
\medskip
\noindent
Saxena, R. K.: 1960, {\sl Proc. Nat. Acad. Sci.
India} {\bf 26}, 400-413\par
\medskip
\noindent
Smith, M.S., Kawano, L.H., Malaney, R.A.: 1993,\par 
Astrophys. J. Suppl. {\bf 85}, 219\par
\medskip
\noindent
Tjon, J., Wu, T.T.: 1979, Phys. Rev. {\bf A19}, 883\par
\medskip
\noindent
Wolfram, S.: 1991, Mathematica - A System for Doing Mathematics\par
by Computer, Addison-Wesley Publishing Company, Inc., Redwood\par
City, California
\end{document}